\documentclass[12pt]{article}
\usepackage{amssymb}
\usepackage{amsfonts}
\usepackage{amsmath}
\usepackage{eufrak}
\usepackage{euscript}
\usepackage{mathrsfs}
\usepackage{verbatim}
\usepackage{fancyhdr}

\def\scr{\EuScript}

\def\orden{{\rm ord}}

\def\expp{{\rm exp}}

\newcommand{\kk}{{\bf K}}
\newcommand{\KK}{{\bf K}}
\newcommand{\QQ}{{\mathbb Q}}
\newcommand{\ZZ}{{\mathbb Z}}
\newcommand{\NN}{{\mathbb N}}

\newcommand{\extL}{{\prec_{\delta}^h}}

\newcommand{\cF}{{\cal F}}
\newcommand{\cG}{{\cal G}}

\newcommand{\cN}{{\cal N}}
\newcommand{\calN}{{\cal N}}

\newcommand{\modul}[1]{\left| #1 \right|}

\newcommand\Z{{\mathbb Z}}

\newcommand\C{{\mathbb C}}
\newcommand\D{{\scr D}}
\newcommand\N{{\mathbb N}}
\newcommand\K{{\bf K}}

\newcommand\expo{\mathop{\rm exp}\nolimits}
\newcommand\Expo{\mathop{\rm Exp}\nolimits}
\newcommand\gr{\mathop{\rm gr}\nolimits}

\newcommand\ord{\prec}
\newcommand\ux{\underline{x}}
\newcommand\uxi{\underline{\xi}}
\newcommand\uchi{\underline{\chi}}
\newcommand\uD{\underline{D}}

\newcounter{numero}[section]

\renewcommand{\thenumero}{\arabic{section}.\arabic{numero}}
\newenvironment{teorema}{\dimen255=\parindent \parindent=0in \medskip
\refstepcounter{numero}{\sc Theorem \thenumero.--}\
\em}{\parindent=\dimen255\par\vspace{1ex}}

\newenvironment{lema}{\dimen255=\parindent \parindent=0in \medskip
\refstepcounter{numero} {\sc Lemma \thenumero.--}\
\em}{\parindent=\dimen255\par\vspace{1ex}}

\newenvironment{definicion}{\dimen255=\parindent \parindent=0in \medskip
\refstepcounter{numero} {\sc Definition \thenumero.--}\
\em}{\parindent=\dimen255\par\vspace{2ex}}

\newenvironment{proposicion}{\dimen255=\parindent \parindent=0in \medskip
\refstepcounter{numero} {\sc Proposition \thenumero.--}\
\em}{\parindent=\dimen255\par\vspace{1ex}}

\newenvironment{nota}{\dimen255=\parindent \parindent=0in \medskip
\refstepcounter{numero}{\sc Remark \thenumero.--}\
\em}{\parindent=\dimen255\par\vspace{2ex}}

\newenvironment{ejemplo}{\dimen255=\parindent \parindent=0in \medskip
\refstepcounter{numero}{\sc Example \thenumero.--}\
\em}{\parindent=\dimen255\par\vspace{2ex}}

\newenvironment{prueba}{\sc Proof:\ \rm}{\hfill $\Box$\par\vspace{2ex}}

\title{Homogenising differential operators\thanks{Supported by DGICYT
PB94-1435.}}
\author{F.J. Castro-Jim\'{e}nez and L. Narv\'{a}ez-Macarro}
\date{June, 1997}
\begin{document}

%
%
%
%
%

\maketitle

\section*{Introduction}

In \cite{laz-83} D. Lazard has used homogenisation of polynomials
to compute the initial ideal $\gr(J)$ of an ideal $J$ generated by
polynomials. In this paper we introduce a homogenisation process
of linear differential operators and we consider ``admissible''
filtrations on the Weyl algebra, generalising $L$-filtrations
\cite{lau_87}. Using an idea similar to Lazard's, we compute
generators of a graded ideal $\gr(I)$ with respect to such
filtrations. As is proved in \cite{a-c-g}, this is a key step to
compute the slopes of a $\cal D$-module.

The Weyl algebra $A_n(\K)$ of order $n$ over a field $\K$ is the
central $\K$-algebra generated by elements $x_i, D_i$,
$i=1,\dots,n$, with relations $[x_i,x_j]=[D_i,D_j]=0, [D_i,x_j]=
\delta_{ij}$. It is naturally filtered by the Bernstein filtration
associated to the total degree in the $x_i$'s and $D_i$'s. Now
consider the {\bf graded} $\K$-algebra $B$, generated by $x_i,
D_i$, $i=1,\dots,n$, and $t$ with homogeneous relations
$[x_i,t]=[D_i,t]=[x_i,x_j]=[D_i,D_j]=0$, $[D_i,x_j]=
\delta_{ij}t^2$. Notice that $A_n(\K)$ is the quotient of $B$ by
the two-sided ideal, generated by the central element $t-1$. In
fact, this algebra coincides with the Rees algebra associated to
the Bernstein filtration of $A_n(\K)$. The homogenisation of an
element in $A_n(\K)$ will be an element in $B$. The homogenisation
process for differential operators we present here has the same
formal properties as the usual homogenisation of commutative
polynomials, and simplifies, considerably, the one studied in
\cite{a-c-g}. We establish the validity of the division theorem
and Buchberger's algorithm to compute standard bases for the
algebra $B$ as in \cite{cas_86}.

Sections 1 and 2 are devoted to the notions of admissible
filtration and $\delta$-standard basis. Section 3 deals with the
main purpose of this paper: homogenisation of differential
operators and effective computation of $\delta$-standard bases.

We would like to thank L\^e D\~ung Tr\'ang for valuable
suggestions.

\section{Admissible filtrations on the Weyl algebra}
\label{filtraciones} Let $\K$ be a field. Let $A_n(\K )$ denote
the Weyl algebra of order $n\geq 1$, i.e.
\begin{eqnarray*}
&A_n(\K) = \K[\ux][\uD]=\K[x_1,\dots,x_n][D_1,\dots,D_n],&\\
&[x_i,x_j]=[D_i,D_j]=0, [D_i,x_j]= \delta_{ij}.&
\end{eqnarray*}

Given a non-zero element $$P= \sum_{\alpha,\beta\in\N^n}
a_{\alpha,\beta}\ux^\alpha \uD^\beta\in A_n(\kk),$$ we denote by
${\cal N}(P)$ its {\em Newton diagram}: $${\cal N}(P) =
\{(\alpha,\beta)\in\N^{2n}\ |\ a_{\alpha,\beta}\neq 0\}. $$

\begin{definicion}
Let $\kk$ be a field. An {\em order function} on $A_n(\kk)$ is a
mapping $\delta: A_n(\kk) \rightarrow \ZZ\cup \{-\infty\}$ such
that:
\begin{enumerate}
\item $\delta(c) = 0$ if $c\in \kk$, $c\neq 0$. \item
$\delta(P)=-\infty$ if and only if $P=0$.
\item $\delta(P+Q)\leq \max\{\delta(P),\delta(Q)\}$. \item $\delta(PQ) =
\delta(P)+\delta(Q)$.
\end{enumerate}
\end{definicion}

\begin{nota}\label{propiedad-1}
If $\delta$ is an order function on $A_n(\kk)$, we have
$\delta(\ux^{\alpha}\uD^{\beta} \ux^{\alpha'}\uD^{\beta'}) =
\delta(\ux^{\alpha+\alpha'}\uD^{\beta+\beta'}).$
\end{nota}

\begin{definicion}
An order function $\delta$ on $A_n(\kk)$ is called {\em
admissible} if, for all non-zero $P\in A_n(\kk)$, we have
$\delta(P) = \max\{\delta(\ux^\alpha\uD^\beta) \,\vert\,
(\alpha,\beta)\in \cN(P)\}$.
\end{definicion}

\begin{proposicion}
Let $\delta: A_n(\kk) \rightarrow \ZZ\cup \{-\infty\}$ be an
admissible order function. Then the family of $\kk$--vector spaces
$$G_\delta^k(A_n(\kk)) = \{P\in A_n(\kk)\ \vert \ \delta(P)\leq
k\}$$ for  $k\in \ZZ$, is an increassing exhaustive separated
filtration of $A_n(\kk)$.
\end{proposicion}

\begin{definicion}
The filtration $G^\bullet_\delta$ will be called the associated
filtration to the admissible order function $\delta$, or the
$\delta$-filtration. A filtration on $A_n(\kk)$ associated to an
admissible order function will be called an {\em admissible
filtration}.
\end{definicion}

\begin{proposicion}
Let $\delta : A_n(\kk) \rightarrow \ZZ\cup \{-\infty\}$ be an
admissible order function. Then the mapping $\Lambda_\delta
:\NN^{2n} \rightarrow \ZZ$ defined  by
$\Lambda_\delta(\alpha,\beta) = \delta(\ux^\alpha\uD^\beta)$ is
the restriction to $\NN^{2n}$ of a unique linear form on
$\QQ^{2n}$, which we will still denote by $\Lambda_\delta$, with
integer coefficients $(p_1,\ldots,p_n,q_1,\ldots,q_n)$ satisfying
$p_i+q_i\geq 0$ for $1\leq i \leq n$.
\end{proposicion}

\begin{prueba} By \ref{propiedad-1},
$\Lambda_\delta(\alpha+\alpha',\beta+\beta') =
\Lambda_\delta(\alpha,\beta) +\Lambda_\delta(\alpha',\beta')$ for
all $\alpha,\alpha',\beta,\beta'
\in
\NN^{n}.$ So there exists a unique linear form $\Lambda_\delta :
\QQ^{2n} \rightarrow \QQ$, with integer coefficients
$(p_1,\ldots,p_n,q_1,\ldots,q_n)$, such that
$\Lambda_\delta(\alpha,\beta) = \delta(\ux^\alpha\uD^\beta)$ for
all $(\alpha,\beta)\in\NN^{2n}$. We have $q_i+p_i = \delta(D_i
x_i) = \delta(x_i D_i +1) = \max\{\delta(x_iD_i),\delta(1)\} =
\max\{p_i+q_i,0\}$ for all $i=1,\dots, n$.
\end{prueba}

\begin{proposicion} \label{0.6}
Let $\Lambda: \QQ^{2n} \rightarrow \QQ$ be a linear form with
integer coefficients $(p_1,\ldots,p_n,q_1,\ldots,q_n)$ satisfying
$p_i+q_i\geq 0$ for all $i=1,\dots, n$. Then there exists a unique
admissible order function $\delta_\Lambda: A_n(\kk) \rightarrow
\ZZ\cup \{-\infty\}$ such that
$\delta_\Lambda(\ux^\alpha\uD^\beta) = \Lambda(\alpha,\beta)$ for
all $(\alpha,\beta)\in\NN^{2n}$.
\end{proposicion}

\begin{prueba}
Let us define  $\delta_\Lambda: A_n(\kk) \rightarrow \ZZ\cup
\{-\infty\}$ by $\delta_\Lambda(0)=-\infty$ and $\delta_\Lambda(P)
= \max\{\Lambda(\alpha,\beta)\ \vert \ (\alpha,\beta)\in \cN(P)\}$
for all non-zero $P\in A_n(\kk)$. Then we have:
\begin{enumerate}
\item $\delta_\Lambda(c) = \Lambda(\underline{0},\underline{0}) = 0$
for all  $c\in \kk$, $c\neq 0$.
\item $\delta_\Lambda(P+Q) = \max \Lambda(\cN(P+Q)) \leq
\max \Lambda(\cN(P)\cup \cN(Q)) = \max (\Lambda(\cN(P)) \cup
\Lambda(\cN(Q))) = \max \{\max \Lambda(\cN(P)), \max
\Lambda(\cN(Q))\} = \max \{\delta_\Lambda(P),
\delta_\Lambda(Q)\}.$
\item For all $i$, $1\leq i \leq n$, we have $\delta_\Lambda(D_ix_i)=
\delta_\Lambda(x_iD_i+1) =
\max\{p_i+q_i,0\}=p_i+q_i=\delta_\Lambda(D_i)+\delta_\Lambda(x_i)$.
\end{enumerate}
Last property implies that $\delta_\Lambda(PQ) =
\delta_\Lambda(P)+\delta_\Lambda(Q)$.
\end{prueba}

\smallskip
Admissible filtrations cover a wide class of filtrations on the
Weyl algebra, as it is showed in the next example.

\begin{ejemplo}
\mbox{ }\\ 1.- The filtration by the order of the differential
operators is the associated filtration to the  admissible order
function $\delta_\Lambda$ where $\Lambda(\alpha,\beta) = |\beta|=
\beta_1+\cdots +\beta_n$. \\ 2.- The $V$-filtration of
Malgrange-Kashiwara with respect to the hypersurface $x_n=0$ is
the associated filtration to the admissible order function
$\delta_\Lambda$ where $\Lambda(\alpha,\beta) = \beta_n-\alpha_n$.
\\ 3.- Let  $L:\QQ^2 \rightarrow \QQ$ be a linear form with
integer coefficients, $L(a,b) = ra+sb$, with $r\geq 0, s\geq 0$
and let us denote by $F_L^\bullet$ the corresponding
$L$-filtration (see \cite{a-c-g}). The $L$-filtration is the
associated filtration to the admissible order function
$\delta_\Lambda$ where $\Lambda(\alpha,\beta) = -s\alpha_n +
r\beta_1+\cdots+r\beta_{n-1} + (r+s)\beta_n$. \\ 4.- For each $k$,
$1\leq k\leq n$, let consider $r, r_1,\ldots,r_k,s_1,\ldots,s_k$
non-negative integers. Let us denote by $\Lambda$ the linear form
on $\QQ^{2n}$ defined by $\Lambda(\alpha,\beta) =
-s_1\alpha_1-\cdots - s_k\alpha_k+(r+s_1)\beta_1 + \cdots +
(r+s_k)\beta_k + r\beta_{k+1}+ \cdots + r\beta_n$. The associated
filtration to the admissible order function $\delta_\Lambda$
coincides with the filtration associated to a multi-filtration
$FV^{\bullet}$ (see \cite{sab_93}).
\end{ejemplo}

\begin{ejemplo} Let $Y$ be the plane curve $x_1^2-x_2^3=0$. The $V$-filtration
on $A_2(\kk)$ associated to $Y$ is not an admissible filtration.
\end{ejemplo}

\begin{nota} \label{grad}
Let $$\Lambda (\alpha,\beta) = p_1\alpha_1+\cdots +p_n\alpha_n +
q_1\beta_1+\cdots +q_n\beta_n $$ be a linear form on $\QQ^{2n}$
with integer coefficients and such that $p_i+q_i\geq 0$ for each
$i=1,\dots,n$. Let $\delta:A_n(\K)\to \Z\cup\{-\infty\}$ the
admissible order function associated to $\Lambda$ (see prop.
\ref{0.6}): $$\delta(P) = \max \{\Lambda(\alpha,\beta)\ |\
(\alpha,\beta)\in\calN(P)\}. $$ If $p_i+q_i>0$ for each
$i=1,\dots,n$, then the graded ring $\gr_\delta(A_n(\K))$
associated to the $\delta$-filtration is commutative. Let us
consider the commutative ring of polynomials
$\K[\chi_1,\dots,\chi_n,\xi_1,\dots,\xi_n]=\K[\uchi,\uxi]$
 and the $\K$-algebra homomorphism
$$\phi: \K[\uchi,\uxi]\longrightarrow \gr_\delta(A_n(\K)), $$ who
sends the $\chi_i$ (resp. the $\xi_i$) to the $\sigma_\delta(x_i)$
(resp. to the $\sigma_\delta(D_i)$), where $\sigma_\delta$ denote
the principal symbol with respect to the $\delta$-filtration. Then
$\phi$ is an isomorphism of graded rings, where $\K[\uchi,\uxi]$
is graded by $$\deg(\uchi^\alpha\uxi^\beta) =
\Lambda(\alpha,\beta).$$

\smallskip
If $p_i+q_i>0$ for each $i=1,\dots,r$, $r<n$, and $p_i+q_i=0$ for
each $i=$ $r+1,\dots,n$, then the graded ring
$\gr_\delta(A_n(\K))$ is non-commutative. Let us consider the
commutative polynomial ring
$R=\K[\chi_1,\dots,\chi_r,\xi_1,\dots,\xi_r]$ and the Weyl algebra
of order $n-r$ over $R$
$$S=R[\chi_{r+1},\dots,\chi_n,\xi_{r+1},\dots,\xi_n]$$ with
relations $$[\chi_i,a]=[\xi_i,a]=0,\quad
[\chi_i,\xi_j]=\delta_{ij} $$ for all $i,j=r+1,\dots,n$ and for
all $a\in R$. The ring $S$ is graded by
$$\deg(\uchi^\alpha\uxi^\beta) = \Lambda(\alpha,\beta) $$ and
there is an isomorphism of graded rings $$\phi:S\longrightarrow
\gr_\delta(A_n(\K))$$ who sends the $\chi_i$ (resp. the $\xi_i$)
to the $\sigma_\delta(x_i)$ (resp. to the $\sigma_\delta(D_i)$).
\end{nota}

\begin{nota} We can also consider admissible filtrations on the ring
$$\D_n=\C\{\ux\}[\uD]$$ of the germs at the origin of linear
differential operators with holomorphic coefficients on $\C^n$. In
this case, admissible order functions $\delta:\D_n\to
\Z\cup\{-\infty\}$ come from linear forms $\Lambda: \QQ^{2n}
\rightarrow \QQ$ with integer coefficients
$(p_1,\ldots,p_n,q_1,\ldots,q_n)$ satisfying $p_i+q_i\geq 0$  and
$p_i\leq 0$ for all $1\leq i\leq n$.
\end{nota}


\newpage

\section{$\delta$-exponents and $\delta$-standard bases in $A_n(\K)$}

Let fix a well monomial ordering $\ord$ in $\N^{2n}$ and let
denote by $\leq$ the usual partial ordering in $\N^{2n}$.

\begin{definicion} For each
admissible order function $\delta:A_n(\K)\to \Z\cup\{-\infty\}$,
we define the following monomial ordering in $\N^{2n}$: $$
(\alpha,\beta)\ord_\delta (\alpha',\beta') \Leftrightarrow \left\{
\begin{array}{l}
\delta(\ux^\alpha\uD^\beta) < \delta(\ux^{\alpha'}\uD^{\beta'})\\
\mbox{or}\ \left\{ \begin{array}{l} \delta(\ux^\alpha\uD^\beta) =
\delta(\ux^{\alpha'}\uD^{\beta'})\\ \mbox{and}\ (\alpha,\beta)\ord
(\alpha',\beta'). \end{array} \right. \end{array} \right. $$
\end{definicion}

\begin{nota} \label{nota-well}
If the order function $\delta$ takes integer negative values, the
ordered set $(\N^{2n},\ord_\delta)$ is not well ordered, but the
restrictions to the level sets of $(\alpha,\beta)$ such that
$\delta(\ux^\alpha\uD^\beta)=c$ are well ordered.
\end{nota}

From now on $\delta:A_n(\K)\to \Z\cup\{-\infty\}$ will denote an
admissible order function.

\begin{definicion} Given a non-zero element $P\in A_n(\K)$,
we define the {\em $\delta$-exponent} of $P$ by $ \displaystyle
\expo_\delta(P) = \max_{\ord_\delta} {\cal N}(P).$ We also denote
by $c_\delta(P)$ the coefficient of the monomial of $P$
corresponding to $\expo_\delta(P)$.
\end{definicion}

\medskip
We have the following classical lemma:

\begin{lema} \label{clasic}
Given two non-zero elements $P,Q$ in $A_n(\K)$ the following
properties hold:
\begin{enumerate}
\item $\expo_\delta(PQ) =\expo_\delta(P) + \expo_\delta(Q)$.
\item If $\expo_\delta(P)\neq \expo_\delta(Q)$ then $\expo_\delta(P+Q) =
\max_{\ord_\delta}\{\expo_\delta(P),\expo_\delta(Q)\}$.
\item If $\expo_\delta(P)=\expo_\delta(Q)$ and if $c_\delta(P)+ c_\delta(Q)\neq
0$ then $\expo_\delta(P+Q) = \expo_\delta(P)$ and
$c_\delta(P+Q)=c_\delta(P)+c_\delta(Q)$.
\item If $\expo_\delta(P)=\expo_\delta(Q)$ and if $c_\delta(P)+c_\delta(Q)=0$
then $\expo_\delta(P+Q)\ord_\delta\expo_\delta(P)$.
\end{enumerate}
\end{lema}

\begin{prueba} 1. We can suppose without loss of generality
$c_\delta(P)=c_\delta(Q)=1$. Let us write $$
\expo_\delta(P)=(\alpha_1,\beta_1), \quad
\expo_\delta(Q)=(\alpha_2,\beta_2) $$ and $$P =
\ux^{\alpha_1}\uD^{\beta_1}+P',\quad Q =
\ux^{\alpha_2}\uD^{\beta_2}+Q'. $$ We have: $$\calN(PQ)\subseteq
\calN(\ux^{\alpha_1}\uD^{\beta_1}\ux^{\alpha_2}\uD^{\beta_2})\cup
\calN(P'\ux^{\alpha_2}\uD^{\beta_2})\cup
\calN(\ux^{\alpha_1}\uD^{\beta_1}Q')\cup \calN(P'Q'). $$ An
element of $\calN(P'\ux^{\alpha_2}\uD^{\beta_2})\cup
\calN(\ux^{\alpha_1}\uD^{\beta_1}Q')\cup \calN(P'Q')$ has the form
$$(\alpha+\gamma-(\beta-\beta'),\beta'+\varepsilon)$$
 with
$\beta'\leq \beta$, $\beta-\beta'\leq\gamma$ and
$(\alpha,\beta)\ord_\delta (\alpha_1,\beta_1),
(\gamma,\varepsilon)\preceq_\delta (\alpha_2,\beta_2)$ or
$(\alpha,\beta)\preceq_\delta (\alpha_1,\beta_1),
(\gamma,\varepsilon)\ord_\delta (\alpha_2,\beta_2)$. By
admissibility
$$\delta(\ux^{\alpha+\gamma-(\beta-\beta')}\uD^{\beta'+\varepsilon})\leq
\delta(\ux^\alpha\uD^\beta\ux^\gamma\uD^\varepsilon) =
\delta(\ux^{\alpha+\gamma}\uD^{\beta+\varepsilon}), $$ and then
$$(\alpha+\gamma-(\beta-\beta'),\beta'+\varepsilon)\preceq_\delta
(\alpha,\beta)+(\gamma,\varepsilon) \ord_\delta
(\alpha_1,\beta_1)+(\alpha_2,\beta_2).$$ Now
$$\calN(\ux^{\alpha_1}\uD^{\beta_1}\ux^{\alpha_2}\uD^{\beta_2})=
\{(\alpha_1+\alpha_2-(\beta_1-\beta'),\beta'+\beta_2)\ |\
\beta'\leq \beta_1,\beta_1-\beta'\leq\alpha_2 \}, $$ and the
monomial $\ux^{\alpha_1+\alpha_2}\uD^{\beta_1+\beta_2}$ can not be
cancelated. So $\expo_\delta (PQ) =
(\alpha_1,\beta_1)+(\alpha_2,\beta_2)$.

\smallskip
The proof of properties 2., 3., 4. is straightforward.
\end{prueba}

\begin{definicion} Given a non-zero left ideal $I$ of $A_n(\K)$,
we define $$\Expo_\delta(I) = \{\expo_\delta (P)\ |\ P\in I, P\neq
0\}.$$
\end{definicion}

\begin{definicion} Given a non-zero left ideal $I$ of $A_n(\K)$,
a {\em $\delta$-standard basis} of $I$ is a family
$P_1,\dots,P_r\in I$ such that $$\Expo_\delta(I) = \bigcup_{i=1}^r
\expo_\delta(P_i)+\N^{2n}.$$
\end{definicion}

\begin{proposicion} \label{bas-st}
Let $I$ be a non-zero left ideal of $A_n(\K)$ and let
$P_1,\dots,P_r$ be a $\delta$-standard basis of $I$. Then
$\sigma_\delta(P_1),\dots,\sigma_\delta(P_r)$ generate the graded
ideal $\gr_\delta(I)$.
\end{proposicion}

\begin{prueba} We follow the proof of lemma 1.3.3 in \cite{a-c-g}.
Let $P$ be a non-zero element of $I$. We define inductively a
family of elements $P^{(s)}$ of $I$ for $s\geq 0$:
\begin{itemize}
\item $P^{(0)}:=P$,
\item $P^{(s+1)} := P^{(s)}-
\frac{c_\delta(P^{(s)})}{c_\delta(P_{i_s})}\ux^{\alpha^s}\uD^{\beta^s}
P_{i_s}$, where $(\alpha^s,\beta^s)$ is an element of $\N^{2n}$
such that $(\alpha^s,\beta^s)+\expo_\delta(P_{i_s}) =
\expo_\delta(P^{(s)})$,
\item $\delta(P^{(s+1)})\leq\delta(P^{(s)})$ and
 $\expo_\delta(P^{(s+1)})\ord_\delta\expo_\delta(P^{(s)})$.
\end{itemize}
By the remark \ref{nota-well}, there is an $s$ such that
$\delta(P^{(s+1)})<\delta(P^{(s)})$. Let $s$ be the smallest
integer having this property. Then $$\sigma_\delta(P) =
\sum_{j=0}^s \sigma_\delta \left(
\frac{c_\delta(P^{(j)})}{c_\delta(P_{i_j})}\ux^{\alpha^j}\uD^{\beta^j}\right)
\sigma_\delta(P_{i_j}). $$
\end{prueba}

\begin{ejemplo} In general, a $\delta$-standard basis of an ideal
$I\subseteq A_n(\K)$ is not a  system of generators of $I$. For
example, take the admissible order function defined by $$\delta
(\ux^\alpha\uD^\beta) = \beta_n-\alpha_n$$ associated to the
$V$-filtration with respect to $x_n=0$, and take $I=A_n(\K)$,
$P=1+x^2_nD_n$.  It is clear that $\sigma_\delta(P)=1$ and then
$P$ is a $\delta$-standard basis of $I$, but obviously $P$ does
not generate $I=A_n(\K)$ ($P$ is not a unit).  \end{ejemplo}

\begin{nota}
Using the isomorphisms $\phi$ of remark \ref{grad}, one can
define, in the obvious way, the Newton diagram $\calN (H)$ for
each non-zero element $H\in\gr_\delta (A_n(\K))$, and so the
$\delta$-exponent $\expo_\delta(H)\in\N^{2n}$, the set
$\Expo_\delta(J)\subseteq \N^{2n}$ for each non-zero left ideal
$J\subseteq \gr_\delta(A_n(\K))$ and the notion of
$\delta$-standard basis for a such $J$.

If $H$ is a non-zero homogeneous element in $\gr_\delta(A_n(\K))$,
then the exponent of $H$ with respect to $\ord$, $\expo_\ord(H)$,
coincides with $\expo_\delta (H)$, and for each non-zero $P\in
A_n(\K)$ we have $$\expo_\delta (P) = \expo_\ord (\sigma_\delta
(P)).$$ In fact, the notions of $\delta$-standard basis and of
$\ord$-standard basis for a non-zero homogeneous left ideal of
$\gr_\delta(A_n(\K))$ coincide.

One can show easily that, for a non-zero left ideal $I\subseteq
A_n(\K)$, a family of elements $P_1,\dots,P_r\in I$ is a
$\delta$-standard basis of $I$ if and only if
$\sigma_\delta(P_1),\dots,\sigma_\delta(P_r)$ is a $\ord$-standard
basis of $\gr_\delta(I)$.

As in \cite{cas_86}, we have a division algorithm in
$\gr_\delta(A_n(\K))$ with respect to the well ordering $\ord$. As
a consequence, a $\ord$-standard basis of $\gr_\delta(I)$ is a
system of generators of this ideal. This precises the proposition
\ref{bas-st}.
\end{nota}

\section{Homog{e}n{i}sation}

In this section we denote by $A_n$ the Weyl algebra $A_n(\K)$. Let
$A_n[t]$  denote the algebra $$A_n[t] =
\K[t,\ux][\uD]=\K[t,x_1,\dots,x_n][D_1,\dots,D_n]$$ with relations
$$[t,x_i]=[t,D_i]=[x_i,x_j]=[D_i,D_j]=0, [D_i,x_j]=
\delta_{ij}t^2.$$ The algebra $A_n[t]$ is graded, the degree of
the monomial $t^k\underline{x}^\alpha\underline{D}^\beta$ being
$k+\modul{\alpha}+\modul{\beta}.$

\begin{lema} The $\KK$--algebra $A_n[t]$ is isomorphic to the Rees
algebra associated to the Bernstein filtration of $A_n$.  The
algebra $\KK[t]$ is central in $A_n[t]$ and the quotient algebra
$A_n[t]/\langle t-1 \rangle$ is isomorphic to $A_n$.  \end{lema}

\begin{prueba} Let $B^\bullet$ be the Bernstein filtration of
$A_n$. We have, for each $m\in \NN$, $$B^m(A_n) =
\{\sum_{\modul{\alpha}+\modul{\beta}\leq m}
p_{\alpha,\beta}\underline{x}^\alpha\underline{D}^\beta\,\vert\,
p_{\alpha,\beta}\in \KK\}.$$ Let $${\cal R}(A_n) =
\bigoplus_{m\geq 0} B^m(A_n)\cdot u^m$$ be the Rees algebra of
$A_n$. We observe that the $\K$-linear map $\phi:A_n[t]
\rightarrow {\cal R}(A_n)$ defined by $$\phi(t)=u,
\phi(x_i)=x_i\cdot u, \phi(D_i)=D_i\cdot u$$ is an isomorphism of
graded algebras.
\end{prueba}

Given  $P=\sum_{\alpha,\beta} p_{\alpha,\beta}
\ux^\alpha\uD^\beta$  in $A_n$ we denote by $\orden^T(P)$ its
total order $$\orden^T(P) = \max\{\modul{\alpha}+\modul{\beta}
\;\vert\; p_{\alpha,\beta}\neq 0\}.$$

\begin{definicion} Let $P= \sum_{\alpha,\beta} p_{\alpha,\beta}
\underline{x}^\alpha\underline{D}^\beta \in A_n.$ Then, the
differential operator $$h(P) = \sum_{\alpha,\beta}
p_{\alpha,\beta} t^{\orden^T(P)-\modul{\alpha}-\modul{\beta}}
\ux^\alpha\uD^\beta \in A_n[t]$$ is called the {\em
homogenisation} of $P$. If
$H=\sum_{k,\alpha,\beta}h_{k,\alpha,\beta}t^k \ux^\alpha\uD^\beta$
is an element of $A_n[t]$ we denote by $H_{|t=1}$ the element of
$A_n$ defined by $H_{|t=1} =
\sum_{k,\alpha,\beta}h_{k,\alpha,\beta}\ux^\alpha\uD^\beta.$
\end{definicion}

\begin{lema}\label{propiedades-de-h}
For $P, Q\in A_n[t]$ we have:
\begin{enumerate}
\item  $h(PQ) = h(P)h(Q)$.
\item There exist $k,l,m\in \NN$ such that $t^kh(P+Q) =
t^lh(P)+t^mh(Q)$.
\end{enumerate}
For any homogeneous element $H\in A_n[t]$ there exists $k\in\NN$
such that $t^kh(H_{|t=1})=H$.
\end{lema}
\begin{prueba}
1. We have $$h(D_ix_i)= h(x_iD_i+1) = x_iD_i + t^2 = D_ix_i =
h(D_i)h(x_i),\quad i=1,\dots,n.$$ From this we obtain easily 1..
To prove 2. let us denote $b=\orden^T(P), c=\orden^T(Q),
d=\orden^T(P+Q)$ and $e=\max\{b,c\}$. We have: $t^{e-d}h(P+Q) =
t^{e-b}h(P) + t^{e-c}h(Q)$. \\ Let $H$ be a non-zero homogeneous
element in $A_n[t]$. Let $k$ be the greatest  integer such that
$t^k$ divides $H$. There exists a homogeneous element $G\in
A_n[t]$ such that $H=t^kG$ and such that the degree of $G$ is
equal to $\orden^T(G_{|t=1})$. We have $H_{|t=1} = G_{|t=1}$ and
$t^kh(G_{|t=1}) = t^kG=H$.
\end{prueba}

Let fix a well monomial ordering $\ord$ in $\N^{2n}$ and an
admissible order function $\delta:A_n\to \Z\cup\{-\infty\}$.

We consider on ${\N}^{2n+1}$ the following total ordering, denoted
by $\extL$, which is a well monomial ordering: $$
\begin{array}{c} (k,\alpha,\beta) {\extL} (k',\alpha',\beta')
\Longleftrightarrow\left\{ \begin{array}{l}
k+\modul{\alpha}+\modul{\beta} < k'+\modul{\alpha'}+\modul{\beta'}
\\ \mbox{or } \left\{ \begin{array}{l}
k+\modul{\alpha}+\modul{\beta} = k'+\modul{\alpha'}+\modul{\beta'}
\mbox{and } \\ (\alpha,\beta) \ord_\delta (\alpha',\beta')
\end{array}\right.  \end{array}\right.  \end{array}$$

\begin{definicion} Let $H=\sum_{k,\alpha,\beta} h_{k,\alpha,\beta}
t^k\underline{x}^\alpha\underline{D}^\beta \in A_n[t]$.  As in
\ref{filtraciones} we denote by $\cN(H)$ the Newton diagram of
$H$: $$\cN(H) = \{(k,\alpha,\beta)\in \NN^{2n+1}\,\vert\,
h_{k,\alpha,\beta}\neq 0\}.$$
\end{definicion}

\begin{definicion}
Given a non-zero element $H\in A_n[t]$ we define the {\em
$\delta$-exponent} of $H$ by $\expo_\delta(H) = \max_{\extL}
\cN(H)$. We also denote by $c_{\delta}(H)$ the coefficient of the
monomial of $H$ corresponding to $\expo_{\delta}(H).$ We write
$\expo(H)$ and  $c(H)$ when no confusion is possible. If $J$ is a
non-zero left ideal of $A_n[t]$ we denote by $\Expo_\delta(J)$ the
set $\{\expo_\delta(H)\ | \ H\in J, H \neq 0\}$.
\end{definicion}

\begin{lema}\label{propiedades-de-exph} Properties 1-4 from lemma \ref{clasic}
hold for $\delta$-exponents of elements in $A_n[t]$. Furthermore,
if $P\in A_n$ then $\pi(\expo_\delta(h(P))) = \expo_\delta(P)$,
where $\pi : \NN^{2n+1}=\NN\times\NN^{2n} \rightarrow \NN^{2n}$ is
the natural projection, and more generally, if $H\in A_n[t]$ is
homogeneous then $\pi(\expo_\delta(H)) =
\pi(\expo_\delta(h(H_{|t=1}))) = \expo_\delta(H_{|t=1})$.
\end{lema}

\begin{prueba} The proof of the first part is similar to that of lemma
\ref{clasic}. Last property follows from lemma
\ref{propiedades-de-h}.
\end{prueba}

\begin{teorema} Let $(P_1,\ldots,P_r)$ be in $A_n[t]^r$.
Let us denote by \begin{eqnarray*} \Delta_1 &:=& \expo (P_1) +
\N^{2n+1}\\ \Delta_i &:=& (\expo (P_i) + \N^{2n+1})\setminus
\bigcup_{j=1}^{i-1} \Delta_j, i=2,\dots ,r\\ \overline{\Delta} &
:= & \N^{2n+1} \setminus \bigcup_{i=1}^r\Delta_i =
\N^{2n+1}\setminus \bigcup_{i=1}^r(\expo (P_i)+\N^{2n+1}).
\end{eqnarray*} Then, for any $H\in A_n[t]$ there exists a unique
element $(Q_1,\ldots,Q_r,R)$ in $A_n[t]^{r+1}$ such that:
\begin{enumerate}
\item $H = Q_1 P_1 + \cdots + Q_r P_r + R.$
\item $\expo(P_i) + \cN(Q_i) \subseteq \Delta_i$
for  $1\leq i \leq r.$
\item $\cN(R) \subseteq \overline{\Delta}.$
\end{enumerate} \end{teorema}

\begin{prueba} The proof is the same as in \cite{cas_86}
since $\extL$ is a well monomial ordering.
\end{prueba}

\begin{lema}\label{proyeccion}
Let $I$ be a non-zero left ideal of $A_n$. We denote by $h(I)$ the
homogenized ideal of $I$, i.e. $h(I)$ is the homogeneous left
ideal of $A_n[t]$ generated by the set $\{h(P)\,\vert\, P\in I\}$.
Then: \begin{enumerate}
\item $\pi(\Expo_\delta(h(I))) = \Expo_\delta(I)$.
\item Let $\{P_1,\ldots,P_m\}$ be a system of generators of $I$. Let
$\widetilde{I}$ be the left ideal of $A_n[t]$ generated by
$\{h(P_1),\ldots,h(P_m)\}$. Then $\pi(\Expo_\delta(\widetilde{I}))
= \Expo_\delta(I)$.
\end{enumerate}
\end{lema}

\begin{prueba}
1. Let $P$ be a non-zero element of $I$. Then the equality
$\expo_\delta(P)
=
\pi(\expo_\delta(h(P)))$ shows that $\Expo_\delta(I) \subseteq
\pi(\Expo_\delta(h(I))).$ Let $H$ a non-zero element of the
homogeneous ideal $h(I)$. We can suppose $H$ homogeneous. There
exist $B_1,\ldots,B_m\in A_n[t]$ and $P_1,\ldots,P_m\in I$ such
that $H= \sum_i B_ih(P_i)$. Then, $H_{|t=1} = \sum_i B_{i|t=1}
P_i$ belong to the ideal $I$. The inclusion
$\pi(\Expo_\delta(h(I)))\subseteq \Expo_\delta(I)$ follows from
\ref{propiedades-de-exph}.\\ 2. Write $P=\sum_iC_iP_i \in I$,
where $C_i\in A_n$. From lemma \ref{propiedades-de-h} there exists
$k\in \NN$ such that $t^kh(P)\in \widetilde{I}$. The equality
$\expo_\delta(P) = \pi(\expo_\delta(t^kh(P)))$ shows that
$\Expo_\delta(I) \subseteq \Expo_\delta(\widetilde{I})$. Finally,
$\widetilde{I} \subseteq h(I)$ implies
$\Expo_\delta(\widetilde{I}) \subseteq \Expo_\delta(h(I)) =
\Expo_\delta(I)$.
\end{prueba}

Let $H_1, H_2$ be elements of $A_n[t]$. Let us denote $\expp(H_i)
= (k_i,\alpha_i,\beta_i)$ and $(k,\alpha,\beta) =
l.c.m.\{(k_1,\alpha_1,\beta_1),(k_2,\alpha_2,\beta_2)\}.$ There
exists $(l_i,\gamma_i,\delta_i)$, for $i=1,2$, such that $
(k,\alpha,\beta) = (k_1,\alpha_1,\beta_1) +
(l_1,\gamma_1,\delta_1) = (k_2,\alpha_2,\beta_2) +
(l_2,\gamma_2,\delta_2).$

\begin{definicion}
The operator $$S(H_1, H_2)=c(H_2)
t^{l_1}x^{\gamma_1}\underline{D}^{\delta_1}H_1 -
c(H_1)t^{l_2}x^{\gamma_2}\underline{D}^{\delta_2} H_2$$  is called
the {\em semisyzygy} relative to $(H_1, H_2)$.
\end{definicion}

\begin{teorema}\label{sicigias}
Let $\cF=\{P_1,\ldots,P_r\}$ be a system of generators of a left
ideal $J$ of  $A_n[t]$ such that, for any $(i,j),$ the remainder
of the division of $S(P_i,P_j)$ by $(P_1,\ldots,P_r)$ is equal to
zero. Then $\cF$ is a $\delta$-standard basis of $J$.
\end{teorema}

\begin{prueba}
This theorem is analogous to Buchberger's criterium for
polynomials \cite{buch_70}. For example, the proof of
\cite[Th.~3.3, Chap.~1]{leje84-85} can be formally adapted to our
case.
\end{prueba}

The preceeding theorem gives an algorithm in order to calculate a
$\delta$-standard basis of an ideal $I$ of $A_n(\K)$ starting from
a system of generators: let $\cF=\{P_1,\ldots,P_r\}$ be a system
of generators of this ideal. We can calculate a $\delta$-standard
basis, say $\cG=\{G_1,\ldots,G_s\}$, of the ideal
$J=\widetilde{I}$ of $A_n[t]$, generated by
$\{h(P_1),\ldots,h(P_r)\}.$ From \ref{proyeccion} we have $
\pi(\Expo_\delta(\widetilde{I})) = \Expo_\delta(I) $ and then, by
\ref{propiedades-de-exph}, $\{G_{1|t=1},\ldots,G_{s|t=1}\}$ is a
$\delta$-standard basis of $I$. Finally, by proposition
\ref{bas-st},
$\{\sigma_\delta(G_{1|t=1}),\ldots,\sigma_\delta(G_{s|t=1})\}$ is
a system of generators of $\gr_\delta(I)$.


\end{document}